\documentclass{amsart}

\usepackage{graphicx}

\usepackage{amsmath,amstext,amssymb,latexsym,lscape}
\usepackage{xcolor}

\theoremstyle{definition}

\theoremstyle{remark}

\numberwithin{equation}{section}

\begin{document}
	
	\title{Energy functional for Lagrangian tori in $\mathbb{C}P^2$}
	
	\author{Hui Ma}
	\address{Department of Mathematical Sciences, Tsinghua University, Beijing 100084, P. R. China}
	\email{hma@math.tsinghua.edu.cn}
	\thanks{The first author was supported in part by NSFC Grant No.~11271213.}

	\author{Andrey E. Mironov}
	\address{Sobolev Institute of Mathematics SB RAS, Pr. acad. Koptyuga 4, 630090, Novosibirsk, Russia and Novosibirsk State University,
   Pirogova Str. 2, 630090, Novosibirsk, Russia}
	\email{mironov@math.nsc.ru}
	\thanks{The second author was supported by the Russian Foundation for Basic Research (grant 16-51-55012) and
  by a grant from Dmitri Zimin's îDynastyî foundation.}
	
	\author{Dafeng Zuo}
	\address{School of Mathematical Science, University of Science and Technology of China, Hefei 230026, P.R. China and Wu Wen-Tsun Key Laboratory of Mathematics, USTC, Chinese Academy of Sciences}
	\email{dfzuo@ustc.edu.cn}
	\thanks{The third author was supported in part by NSFC (Grant No. ~11671371,11371338) and the Fundamental Research
Funds for the Central Universities.}

	
	\date{\today}
	
	
	
	\begin{abstract}
		In this paper we study Lagrangian tori in ${\mathbb C}P^2$. A two-dimensional periodic Schr\"odinger operator
is associated with every Lagrangian torus in ${\mathbb C}P^2$. We introduce an energy functional for tori as an integral of the potential
of the Schr\"odinger operators, which  has a natural geometrical meaning (see below). We study the energy functional on two families
of Lagrangian tori and propose a conjecture that the minimum of the functional is achieved by the Clifford torus.
We also study deformations of minimal Lagrangian tori. In particular we show that if the deformation preserves a conformal type of the torus,
then it also preserves the area of the torus. Thus it follows that deformations generated by Novikov--Veselov equations preserve the area of minimal Lagrangian tori.
\end{abstract}
	
	\maketitle

\section{Introduction and main results}

In this paper we study Lagrangian tori in ${\mathbb C}P^2$. The paper consists of two parts. In the first part we introduce an energy functional of
Lagrangian torus as an integral
of the potential of Schr\"odinger operator associated with the torus. We study this energy functional on a family of homogeneous tori and 
Hamiltonian-minimal Lagrangian tori constructed in \cite{M3}, respectively.
 We propose a conjecture that the minimum of this functional on the whole set of Lagrangian tori can be achieved by the Clifford torus. Similar conjecture
 was proposed by Montiel and Urbano \cite{MU} for their Willmore functional. In the second part of the paper we study deformations of minimal
 Lagrangian tori. Such deformations are related
 to eigenfunctions of the Laplace--Beltrami operator with eigenvalue $\lambda=6$. Using Novikov--Veselov  hierarchy \cite{NV}
 we propose a method of finding such eigenfunctions.
  We also prove that if the deformation of minimal Lagrangian tori preserves the conformal type of the torus then it preserves the area of the torus.


Let $\Sigma$ be a closed Lagrangian surface immersed in $\mathbb{C}P^2$.
Let $x,y$ denote local conformal coordinates such that the induce metric of $\Sigma$ is given by
\begin{equation}\label{eq1}
 ds^2=2e^{v(x,y)}(dx^2+dy^2).
\end{equation}
Let  $r: U\rightarrow S^5$ be a local horizontal lift of the immersion defined on an open subset $U$ of $\Sigma$, where $S^5\subset {\mathbb C}^3$ is the
unit sphere.
Since $\Sigma$ is Lagrangian and $x,y$ are conformal coordinates we have
\begin{equation}\label{eq2}
\langle r, r\rangle=1, \ \langle r_x, r\rangle=\langle r_y, r\rangle=\langle r_x,r_y\rangle =0,
\end{equation}
\begin{equation}\label{eq3}
 \langle r_x, r_x\rangle =\langle r_y, r_y \rangle=2e^v,
\end{equation}
and hence
$$\tilde{R}=\begin{pmatrix}r\\ \frac{r_x}{|r_x|}\\ \frac{r_y}{|r_y|} \end{pmatrix}\in {\rm U}(3),$$
where $\langle\, ,\, \rangle$ denotes the Hermitian inner product of $\mathbb{C}^3$.
One can define a local function $\beta: U\rightarrow \mathbb{R}$, called the {\emph Lagrangian angle} of $\Sigma$,
by $e^{i\beta(x,y)}=\det \tilde{R}$.
Consequently
$$R=\begin{pmatrix}r\\ \frac{1}{\sqrt{2}} e^{-\frac{v}{2}-i\frac{\beta}{2}} r_x\\\frac{1}{\sqrt{2}} e^{-\frac{v}{2}-i\frac{\beta}{2}} r_y \end{pmatrix}\in {\rm SU}(3).$$
By direct calculations, by
\begin{equation}\label{eq3.1}
R_x=AR, \qquad R_y=BR,
\end{equation}
where $A,B\in {\rm su}(3)$ (see Section 3)
we obtain the following lemma.

\bigskip
\noindent
{\bf Lemma 1.}  \cite{M1} {\it
The vector function $r$ satisfies the Schr\"{o}dinger equation $Lr=0$,
where
$$L=(\partial_x -\frac{i\beta}{2})^2+(\partial_y-\frac{i\beta}{2})^2+V(x,y),$$
with the potential
$$
 V=4e^v+\frac{1}{4} (\beta_x^2+\beta_y^2)+\frac{i}{2}\triangle \beta.
$$}

\bigskip

Let us assume that $\Sigma$ is a torus given by the mapping
$$
 r:{\mathbb R}^2\rightarrow S^5,
$$
where $r$ satisfies \eqref{eq2} and  \eqref{eq3}. Then the potential $V$ is a doublely periodic function with respect to a  lattice of
periods
$\Lambda\subset \mathbb{R}^2$ and $r$ is a Bloch eigenfunction of the Schr\"{o}dinger operator $L$, i.e.,
$$r((x, y)+\lambda_s)=e^{i p_s} r(x,y), \, \ p_s\in \mathbb{R}, \, \ s=1,2,$$
where $\lambda_1, \lambda_2$ is a basis of $\Lambda$.
Let us introduce the energy  of a Lagrangian torus as an integral of $V$ (Remark that
similarly one can define the energy functional
of arbitrary Lagrangian surfaces in $\mathbb{C}P^2$).

\bigskip
\noindent
{\bf Definition 1.} {\it The energy of the Lagrangian torus  $\Sigma\subset\mathbb{C}P^2$ is
$$E(\Sigma)=\frac{1}{2}\int_{\Sigma} V dx\wedge dy.$$}
\bigskip

\noindent
{\bf Remark 1.} {\it For convenience, we take the factor $\frac{1}{2}$ in the above definition.}

\bigskip

\noindent
It turns out that the energy functional has the following geometric meaning.

\bigskip
\noindent
{\bf Lemma 2.} {\it The energy of the Lagrangian torus is
$$E(\Sigma)=\int_{\Sigma} d\sigma +\frac{1}{8} \int_{\Sigma} |H|^2 d\sigma,$$
where $d\sigma=2e^v dx\wedge dy$ is the area element of $\Sigma$, $H$ is the mean curvature vector.}


\bigskip

\noindent
Let $\Sigma_{r_1,r_2,r_3}$ be a homogeneous torus in ${\mathbb C}P^2$, where $r_1,r_2,r_3$ are positive numbers such that $r_1^2+r_2^2+r_3^2=1$.
Any homogeneous torus can be obtained
as the image of the Hopf projection
$\mathcal{H}:S^5\rightarrow {\mathbb C}P^2$
of the 3-torus
$$\{(r_1e^{i\varphi_1},r_2e^{i\varphi_2},r_3e^{i\varphi_3}), \varphi_j\in{\mathbb R}\}\subset S^5.$$
 Every Homogeneous torus is Hamiltonian-minimal Lagrangian in ${\mathbb C}P^2$, i.e.
its area is a critical point under Hamiltonian deformations (see \cite{O}).
Among homogeneous tori there is a minimal torus $\Sigma_{\frac{1}{\sqrt{3}},\frac{1}{\sqrt{3}},\frac{1}{\sqrt{3}}}$, which is called Clifford torus
(we denote it by $\Sigma_{Cl}$).

\bigskip

\noindent
{\bf Proposition 1.} {\it
 The identity
 $$
  E(\Sigma_{r_1,r_2,r_3})=\frac{\pi^2 (1 - r_1^2) (1 - r_2^2) (1-r_3^2)}{2 r_1 r_2 r_3}
 $$
 holds.
 Among Lagrangian homogeneous tori in $\mathbb{C}P^2$, the Clifford torus attains the minimum of the energy functional
 $$E(\Sigma_{Cl})=Area(\Sigma_{Cl})=\frac{4 \pi^2}{3 \sqrt{3}}.$$
}

Let us consider another family of Hamiltonian-minimal Lagrangian tori in $\mathbb{C}P^2$ (see \cite{M3}).
Let $\Sigma_{m,n,k}\subset {\mathbb C}P^2$ denote a surface given as the image of the surface
$$
 \{u_1e^{2\pi imy},u_2e^{2\pi iny},u_3e^{2\pi iky}\}\subset S^5
$$
under the Hopf projection, where $(u_1,u_2,u_3)\in {\mathbb R}^3$ such that
$$
 u_1^2+u_2^2+u_3^2=1,
$$
$$
 mu_1^2+nu_2^2+ku_3^2=0,
$$
with integers $m\geq n>0$ and $k<0$. Such surface $\Sigma_{m,n,k}$ is an (immersed or embedded) Hamiltonian-minimal Lagrangian torus or Klein bottle. The topology of
$\Sigma_{m,n,k}$ depends on whether the involution
$$
 (u_1,u_2,u_3)\rightarrow (u_1\cos(m\pi),u_2\cos(n\pi),u_3\cos(k\pi))
$$
preserves an orientation of the surface $mu_1^2+nu_2^2+ku_3^2=0$.
We have

\bigskip

\noindent
{\bf Proposition 2.} {\it The energy of $\Sigma_{m,n,k}$ is greater than the energy of the Clifford torus
$$E(\Sigma_{m,n,k}) > E(\Sigma_{Cl}).$$}
\bigskip

Propositions 1 and 2 allow us to propose the following conjecture.
\bigskip

\noindent
{\bf Conjecture. }
{\it The Clifford torus attains the minimum of the energy functional among all Lagrangian tori in $\mathbb{C}P^2$.}
\bigskip

\noindent

In the second part of the paper we study deformations of minimal Lagrangian tori in ${\mathbb C}P^2$.

\bigskip
\noindent
{\bf Theorem 1.} {\it
Let $\Sigma_0$ be a closed minimal Lagrangian torus in $\mathbb{C}P^2$
and $\Sigma_t$ be a smooth deformation of $\Sigma_0$ preserving the conformal type of $\Sigma_0$ such that $\Sigma_t$ is still minimal Lagrangian.
 Then the area of $\Sigma_t$ is preserved
 $$\mathrm{Area}(\Sigma_t)=\mathrm{Area}(\Sigma_0).$$
}
\bigskip

Important examples of such deformations are deformations generated by Novikov--Veselov hierarchy.
Let $\Sigma\subset{\mathbb C}P^2$ be a minimal Lagrangian torus. On $\Sigma$ there are coordinates $x,y$ such that the induced metric has the form  $ds^2=2e^{v(x,y)}(dx^2+dy^2),$ where $v$ satisfies
the Tziz\'{e}ica equation
\begin{equation}\label{eq4}
 \Delta v=4(e^{-2v}-e^{v}).
\end{equation}
Let $r:{\mathbb R}^2\rightarrow S^5$ be a lift map for $\Sigma$.

\bigskip
\noindent
{\bf Theorem 2.} \cite{M1} {\it There is a mapping $\tilde{r}(t):{\mathbb R}^2\rightarrow S^5$, $t=(t_2,t_3,\dots),$ $\tilde{r}(0)=r$,
defining a minimal Lagrangian torus $\Sigma_t\subset{\mathbb C}P^2$ such that $\Sigma_0=\Sigma$.
The map $\tilde{r}$ satisfies the equations
$$
 L\tilde{r}=\Delta\tilde{r}+4 e^{\tilde{v}}\tilde{r},
$$
$$
 \partial_{t_n}\tilde{r}=A_n\tilde{r},
$$
where $A_n$ are operators of order $(2n+1)$ in $x,y$ and $ds^2=2e^{\tilde{v}(x,y,t)}(dx^2+dy^2)$ is the induced metric on $\Sigma_t$.
The deformation $\tilde{r}(t)$ preserves conformal type of the torus and the spectral curve of the Schr\"{o}dinger operator $L$.
The function $\tilde{V}=4 e^{\tilde{v}}$ with $\tilde{v}(0)=v$
satisfies  the Novikov--Veselov hierarchy
$$
 \frac{\partial{L}}{\partial_{t_n}}=[L,A_n]+B_nL.
$$}

\bigskip
Thus Theorems 1 and 2 leads to the following corollary.

\bigskip
\noindent
{\bf Corollary 1.} {\it Deformations of minimal Lagrangian tori given by Novikov--Veselov hierarchy (see Theorem 2) preserve the area of tori.}
\bigskip

It would be interesting to use Novikov--Veselov equations to construct deformations preserving the energy of arbitrary Lagrangian torus in
${\mathbb C}P^2$. Similar construction is possible for tori in ${\mathbb R}^3$ if one use modified Novikov--Veselov equations (see \cite{T1}, \cite{T2}).

Deformations of minimal Lagrangian tori are related to eigenfunctions of the Laplace--Beltrami operator with eigenvalue 6.
Generally to find explicitly eigenfunctions of the Laplace--Beltrami operator is a hard problem (see \cite{P}).
The preceding theorem gives a method to find such eigenfunctions for Laplace--Beltrami operator of minimal Lagrangian tori in ${\mathbb C}P^2$.

\bigskip
\noindent
{\bf Theorem 3.} {\it Let $ds^2=2e^{v(x,y)}(dx^2+dy^2)$ be a metric on a surface $\Sigma$ and $v$ satisfies the Tziz\'{e}ica equation.
Then functions
$$
 s_1=v_x^2-v_y^2+v_{xx}-v_{yy},
$$
$$
 s_2=v_xv_y+v_{xy}
$$
are eigenfunctions of the Laplace--Beltrami operator with eigenvalue 6.}

Other eigenfunctions can be found with the help of explicit calculations related to Novikov--Veselov equations.

\section{Energy functional}

In this section we prove Lemma 1,  Propositions 1 and 2.

\subsection{Proof of Lemma 2}

Recall that the mean curvature vector field on Lagrangian surface $\Sigma\subset\mathbb{C}P^2$ can be expressed in terms of Lagrangian angle
$\beta$ (see for example \cite{M3})
$$H=J\mathrm{grad} \beta,$$
where $J$ is the standard complex structure on $\mathbb{C}P^2$. With respect to the induced metric $ds^2=2e^v (dx^2+dy^2)$, 
$$|H|^2=|J\mathrm{grad} \beta|^2=|\mathrm{grad}\beta|^2=\frac{e^{-v}}{2}(\beta_x^2+\beta_y^2).$$
Hence
\begin{eqnarray*}
E(\Sigma)&=&\int_{\Sigma} 2e^v dx\wedge dy + \frac{1}{8} \int_{\Sigma} 2e^v |H|^2 dx\wedge dy + \frac{i}{4}\int_{\Sigma} \triangle \beta dx\wedge dy\\
&=&\int_{\Sigma} d\sigma +\frac{1}{8}\int_{\Sigma} |H|^2 d\sigma.
\end{eqnarray*}


\subsection{Proof of Proposition 1}
Let $\Sigma_{r_1,r_2,r_3}\subset \mathbb{C}P^2$ be a homogeneous torus.
Then after choosing of appropriate coordinates and taking automorphisms of $\mathbb{C}P^2$, the horizontal
lift $r: \mathbb{R}^2 \rightarrow S^5$ of $\Sigma_{r_1,r_2,r_3}$ can be given by
$$
 r(x,y)=(r_1 e^{2 \pi i x} , r_2 e^{2 \pi i (a_1 x + b_1 y)}, r_3 e^{2 \pi i (a_2 x + b_2 y)}),
$$
where $r_1^2+r_2^2+r_3^2=1$.  It follows from \eqref{eq2} and \eqref{eq3} that
$$
 a_1=a_2=-\frac{r_1^2}{r_2^2+r_3^2}, \quad b_1=\frac{r_1 r_3}{r_2(r_2^2+r_3^2)},
 \quad  b_2=-\frac{r_1 r_2}{r_3(r_2^2+r_3^2)}.
$$
By direct calculations, one can obtain that the lattice of periods
$\Lambda$ for $\mathcal{H}\circ r$,
is $\Lambda=\{\mathbb{Z} e_1+\mathbb{Z}  e_2 \}\subset \mathbb{R}^2$, where
$$e_1=(r_2^2+r_3^2,0),  \quad e_2=(r_3^2, \frac{r_2 r_3}{r_1}).$$
Moreover, it follows from
$\langle r_x, r_x\rangle=\langle r_y, r_y\rangle=\frac{4\pi^2 r_1^2}{r_2^2+r_3^2} $
 that the induced metric, the area, the Lagrangian angle and the potential of the Schr\"{o}dinger operator (see Lemma 1)  of $\Sigma_{r_1,r_2,r_3}$ are given by
\begin{eqnarray*}
ds^2&=&\frac{4\pi^2 r_1^2}{r_2^2+r_3^2} (dx^2+dy^2),\\
\int_{\Sigma_{r_1,r_2,r_3}}d\sigma&=& 4\pi^2 r_1 r_2 r_3,\\
 \beta&=& 2\pi \frac{ 1-3r_1^2}{r_2^2+r_3^2} x-2\pi\frac{ r_1(r_2^2-r_3^2)}{r_2 r_3 (r_2^2+r_3^2)}y,\\
 V&=&\frac{\pi^2 (1-r_2^2)(r_1^2+r_2^2)}{r_2^2 r_3^2}.
\end{eqnarray*}
Hence the energy of $\Sigma_{r_1,r_2,r_3}$ is given by
\begin{eqnarray*}
 E(\Sigma_{r_1,r_2,r_3})&=&4\pi^2 r_1 r_2 r_3+\pi^2 \frac{(r_1^2 r_2^4 + r_2^2 r_3^2 - 8 r_1^2 r_2^2 r_3^2 +
   9 r_1^4 r_2^2 r_3^2 + r_1^2 r_3^4)}{2 r_1 r_2 r_3 (r_2^2 + r_3^2)}\\
   &=&\frac{\pi^2 (1 - r_1^2) (1 - r_2^2) (1-r_3^2)}{2 r_1 r_2 r_3}.
\end{eqnarray*}

Thus it is easy to find that the minimum of $E(\Sigma_{r_1,r_2,r_3})$ is attained on Clifford torus $\Sigma_{\frac{1}{\sqrt{3}},\frac{1}{\sqrt{3}},\frac{1}{\sqrt{3}}}.$

\subsection{Proof of Proposition 2}

Taking an appropriate parametrization of the curve given by
\begin{eqnarray*}
u_1^2+u_2^2+u_3^2=1,\\
m u_1^2+nu_2^2+k u_3^2=0,
\end{eqnarray*}
where $m\geq n>0$ and $k<0$ are constant integers,
we obtain
$$\psi(u,y)=(u_1 e^{2\pi i m y}, u_2 e^{2\pi i n y}, u_3 e^{2\pi i ky}) \subset S^5,$$
with
$$u_1=\sin x \sqrt{\frac{k}{k-m}}, u_2=\cos x \sqrt{\frac{k}{k-n}}, u_3=\sqrt{\frac{n\cos^2 x}{n-k}+
\frac{m\sin^2 x}{m-k}}.$$
Thus we have a horizontal lift of a surface $\Sigma$ in $\mathbb{C}P^2$ given by
$$r(x,y)=(u_1(x)  e^{2\pi i m y}, u_2(x) e^{2\pi i n y}, u_3(x) e^{2\pi i ky}).$$
Now let us consider the torus given as above, denoted it by $\Sigma_{m,n,k}$.

By straightforward calculations, we obtain the induced metric on $\Sigma_{m,n,k}$ is
$$ds^2=2e^{v_1(x)} dx^2 + 2e^{v_2(x)} dy^2,$$
where
\begin{eqnarray*}
2e^{v_1(x)}&=&-\frac{k(m+n-(m-n)\cos (2x))}{2mn-k(m+n)+k(m-n) \cos (2x)},\\
2e^{v_2(x)}&=&-2k\pi^2 (m+n-(m-n)\cos(2x)).
\end{eqnarray*}

Now we have
$$\tilde{R}=\begin{pmatrix}
r\\
\frac{r_x}{|r_x|}\\
\frac{r_y}{|r_y|}
\end{pmatrix}\in {\rm U}(3).$$
By the definition of the Lagrangian angle $\det \tilde{R}=e^{i\beta}$, we get
$$e^{i\beta}=ie^{2\pi i (m+n+k)y}.$$

The period in $x$ is $e_1=2\pi$ and the period in $y$ is $\frac{1}{p}$, where $p=(m-k,n-k)$ is the biggest common factor of $m$ and $n$. Moreover, the Willmore functional is given by



\begin{eqnarray*}
W(\Sigma_{m,n,k})&=&\int_{\Sigma_{m,n,k}} |H|^2 d\sigma=\int_{[0,2\pi]\times [0,\frac{1}{p}]} |H|^2 2e^{\frac{v_1+v_2}{2}} dxdy\\
&=&\int_{\Sigma_{m,n,k}} (e^{\frac{v_2-v_1}{2}}\beta_x^2 +e^{\frac{v_1-v_2}{2}}\beta_y^2 )dxdy.
\end{eqnarray*}

Thus we have
\begin{eqnarray*}
A(\Sigma_{m,n,k})&=&\int_{\Sigma_{m,n,k}} d\sigma=\frac{1}{p}\int_{0}^{2\pi} \frac{-\sqrt{2} k\pi (m+n-(m-n)\cos 2x)}{\sqrt{2mn-k(m+n)+k(m-n)\cos 2x}} dx,\\
W(\Sigma_{m,n,k})&=&\frac{1}{p}\int_0^{2\pi} \frac{2\sqrt{2} (k+m+n)^2 \pi}{\sqrt{2mn-k(m+n)+k(m-n)\cos 2x}} dx,\\
E(\Sigma_{m,n,k})&=& A(\Sigma_{m,n,k})+\frac{1}{8} W(\Sigma_{m,n,k})\\
&=&\frac{1}{p} \frac{\pi}{2\sqrt{2}}  \int_0^{2\pi} \frac{4k(m-n)\cos 2x +(-k+m+n)^2}{\sqrt{2mn-k(m+n)+k(m-n)\cos 2x}} dx.
\end{eqnarray*}
Remembering $k<0$, we have
$$\frac{4k(m-n)\cos 2x +(-k+m+n)^2}{\sqrt{2mn-k(m+n)+k(m-n)\cos 2x}} \geq
\frac{4k(m-n)\cos 2x +(-k+m+n)^2}{\sqrt{2mn-k(m+n)-k(m-n)}}.$$
Hence
\begin{eqnarray*}
E(\Sigma_{m,n,k})&\geq& \frac{1}{p} \frac{\pi}{2\sqrt{2}}  \int_0^{2\pi} \frac{4k(m-n)\cos 2x +(-k+m+n)^2}{\sqrt{2mn-k(m+n)-k(m-n)}}\\
&=&  \frac{\pi^2}{2p}  \frac{(-k+m+n)^2}{\sqrt{m(-k+n)}}\\
&\geq& \frac{\pi^2}{p} (-k+n+m)=\frac{m+pr}{p} \pi^2\\
&=&(\frac{m}{p}+r)\pi^2 >\frac{4}{3\sqrt{3}} \pi^2,
\end{eqnarray*}
where we use $\frac{4}{3\sqrt{3}}\sim 0.8$, $n-k=pr$ and $r$ is a positive integer.
Thus the proof of Proposition 1.7 is completed.

\section{deformations of Minimal Lagrangian tori}
Now consider a Lagrangian torus defined by the composition of the maps
$$
 r:{\mathbb R}^2\rightarrow S^5
$$
and ${\mathcal H}$. Then $R$ satisfies \eqref{eq3.1} (see \cite{M1}), where
$$A=\begin{pmatrix}
0&\sqrt{2} e^{\frac{v+i\beta}{2}}&0\\
-\sqrt{2} e^{\frac{v-i\beta}{2}} & iFe^{-v}& -\frac{v_y}{2}+i(e^{-v}G+\frac{1}{2} \beta_y)\\
0& \frac{v_y}{2}+i(e^{-v}G+\frac{1}{2}\beta_y)  &-iF e^{-v}
\end{pmatrix}
\in {\rm su}(3),
$$
$$
B=\begin{pmatrix}
0&0&\sqrt{2} e^{\frac{v+i\beta}{2}}\\
0& -i Ge^{-v} & \frac{1}{2} v_x+i(-e^{-v}F+\frac{1}{2}\beta_x)\\
-\sqrt{2} e^{\frac{v-i\beta}{2}} &-\frac{1}{2}v_x+i(-e^{-v} F+\frac{1}{2}\beta_x) &-iG e^{-v}
\end{pmatrix}
\in {\rm su}(3),
$$
with real functions $F$ and $G$ given by
$$F=-\frac{1}{2i} (\langle r_{xy}, r_y\rangle -e^{v}(v_x+i\beta_x))$$ and
$$G=\frac{1}{2i}  ( \langle r_{xy}, r_{x}\rangle -e^v (v_y+i\beta_y)).$$
The compatibility condition
$$A_y-B_x+[A,B]=0$$
 leads to the following equations (see \cite{M1} and also \cite{Ma1}, \cite{Ma2})
$$
2F_x+2 G_y=(\beta_{xx}-\beta_{yy})e^v,
$$
$$
2F_y-2G_x=(\beta_x v_y+\beta_y v_x)e^v,
$$
$$
\triangle v=4(F^2+G^2)e^{-2v}-4e^v-2 (F\beta_x-G\beta_y) e^{-v}.
$$
If the torus is minimal then $\beta=const$ and from the equations it follows that $F$ and $G$ are constants. After appropriate change of
coordinates (a homothety and rotation) one can assume that $F=1$ and $G=0$. Hence $v$ satisfies the Tziz\'{e}ica equation (\ref{eq4}).
Smooth periodic solutions of the Tziz\'{e}ica equation are finite-gap solutions. These solutions were found in \cite{Sh}.
Minimal Lagrangian tori were studied in \cite{Ma3}--\cite{J}.
Assume that we have a deformation $\Sigma_t$ of $\Sigma$, $\Sigma_0=\Sigma,$ given by the mapping
$$
 r(t):{\mathbb R}^2\rightarrow S^5
$$
with the induced metric $ds^2=2e^{v(x,y,t)}(dx^2+dy^2),$ $v(x,y,t)$ satisfies (\ref{eq4}).
We have
 $R_t=T R$, where
\begin{equation}\label{eq:C}
 T=\begin{pmatrix}
i s & a_1+ib_1 & a_2+ib_2\\
-a_1+ib_1& is_1&a_3+ib_3\\
-a_2+ib_2 & -a_3 +ib_3& -i(s+s_1)
\end{pmatrix} \in {\rm su}(3),
\end{equation}
with functions $s, s_1, a_1,a_2,a_3, b_1, b_2$ and $b_3$ depending on $x,y$ and $t$.
From the compatibility conditions $$A_t-T_x+[A,T]=0, \quad B_t-T_y+[B,T]=0,$$
we obtain the identities
$$
b_1=\frac{e^{-\frac{v}{2}}s_x}{2\sqrt{2}}, \quad \quad \quad b_2=\frac{e^{-\frac{v}{2}} s_y}{2\sqrt{2}},
$$
$$
a_3=\frac{e^{-2v}(\sqrt{2}s_y-2e^{\frac{3}{2}v}(a_1 v_y- 2 {a_2}_x))}{4\sqrt{2}},
$$
$$
b_3=\frac{e^{-\frac{3}{2}v}(-8 a_2 - \sqrt{2} e^{\frac{1}{2}v} ( s_x v_y + s_y v_x -2s_{xy}) )}{8\sqrt{2}},
$$
$$
s_1=\frac{e^{-\frac{3}{2}v} (8a_1+\sqrt{2} e^{\frac{1}{2}v} (8e^v s+s_y v_y-s_x v_x+2s_{xx}))}{8\sqrt{2}}
$$
and the following overdetermined system of equations for $s, a_1, a_2$ which determine the deformation:
\begin{equation}\label{eq5}
 v_t+ \frac{1}{2} e^{-2v}(-\sqrt{2}e^{\frac{3}{2}v} a_2 v_y-2\sqrt{2} e^{\frac{3}{2}v} {a_1}_x+s_x)=0,
\end{equation}
\begin{equation}\label{eq6}
 \triangle s+ 12 e^v s=0,
\end{equation}
\begin{equation}\label{eq7}
2 {a_1}_x-2 {a_2}_y -a_1 v_x  + a_2 v_y -\sqrt{2}e^{-\frac{3}{2}v} s_x =0,
\end{equation}
\begin{equation}\label{eq8}
2 {a_1}_y +2 {a_2}_x - a_1 v_y - a_2 v_x   + \sqrt{2} e^{-\frac{3}{2}v} s_y=0,
\end{equation}
\begin{equation}\label{eq9}
 \sqrt{2} e^{\frac{3v}{2}} a_2 v_y+2\sqrt{2} e^{\frac{3v}{2}} {a_1}_x -s_x
+ 2e^{3v} (b_3 v_y+s_x + {s_1}_x)=0,
\end{equation}
\begin{equation}\label{eq10}
 2e^{-v} a_3-{s_1}_y+b_3 v_x=0.
\end{equation}
Now we can give a geometric interpretation of $s$. Let $r_t^{\bot}$ be the normal component of the velocity vector $r_t$. Thus
$$
 r_t^{\bot} = \frac{b_1}{\sqrt{2} e^{\frac{v}{2}}} i r_x +\frac{b_2}{\sqrt{2} e^{\frac{v}{2}}} i r_y
 =\frac{e^{-v}}{4} (s_x i r_x + s_y i r_y)=\mathrm{grad} \frac{s}{2}.
$$
Thus it follows that the deformation we obtained above is a Hamiltonian deformation with Hamiltonian $\frac{s}{2}$ (see \cite{O}).
Moreover (\ref{eq6}) can be rewritten in the following form
$$
 \triangle_{_{LB}}s=6s,
$$
where $\triangle_{_{LB}}$ is the Laplace--Beltrami operator of the metric (\ref{eq1}). Thus the function $s$ is an eigenfunction of the Laplace--Beltrami operator with the eigenvalue $6$.

\begin{proof}[Proof of Theorem 1]

Using (\ref{eq5}) and (\ref{eq7}), we get
\begin{eqnarray*}
v_t&=&\frac{e^{-\frac{v}{2}} (2{a_2}_y +a_2 v_y +2 {a_1}_x +a_1 v_x)}{2\sqrt{2}},\\
s_x&=&-\frac{e^{\frac{3}{2}v} (2 {a_2}_y-a_2 v_y-2{a_1}_x +a_1 v_x)}{\sqrt{2}}.
\end{eqnarray*}
Considering the area form $d\sigma=2e^v dx\wedge dy$, we set
$$\Omega=\partial_ t(2 e^v )dx\wedge dy
=\frac{e^{\frac{v}{2}} (2 {a_2}_y+a_2 v_y+2{a_1}_x+a_1 v_x)} {\sqrt{2}} dx\wedge dy.$$
It turns out that  $\Omega=d\omega$,
where
$$\omega= \sqrt{2} e^{\frac{v}{2}} (a_1dy -a_2 dx).$$
If $\Sigma_t$ is a smooth deformation of $\Sigma_0$ preserving the conformal type of $\Sigma_0$,
then
$$ \frac{d}{dt}  \int_{\Lambda} 2e^v dx\wedge dy =\int_{\Lambda} \Omega =\int_{\Lambda} d\omega =0,$$
where $\Lambda$ is a lattice of periods.
The proof is completed.

\end{proof}

\subsection{Novikov--Veselov hierarchy and deformations of minimal Lagrangian tori}
In this subsection we consider an example of deformation of minimal Lagrangian tori defined by the second Novikov--Veselov
equation. In particular we give an explicit solution of the system (\ref{eq5})--(\ref{eq10}) in
terms of the function $v$ defining the induced metric (\ref{eq1}) of the torus.

Let us recall the Novikov--Veselov hierarchy.
Let $L$ be a Schr\"{o}dinger operator
$$L=\partial_z\partial_{\bar z}+V(z,\bar{z}),$$
and
$$A_{2n+1}= \partial_z^{2n+1} +u_{2n-1}\partial_z^{2n-1}+\cdots+u_1\partial_z+
\partial_{\bar{z}}^{2n+1}+w_{2n-1}\partial_{\bar z}^{2n-1}+\cdots+w_1\partial_{\bar z},$$
where $u_j=u_j(z,\bar{z}), w_j=w_j(z,\bar{z})$.
The operator $\partial_{t_n}-A_{2n+1}$ defines an evolution equation for the eigenfunction $r$ of the Schr\"{o}dinger operator $Lr=0$ by
$$\partial_{t_n} r=A_{2n+1} r. $$
The Novikov--Veselov equations are
$$\frac{\partial L}{\partial t}=[A_{2n+1}, L]+B_{2n-2} L,$$
where $B_{2n-2}$ is a differential operator of order $2n-2$ and
 $V$ and coefficients of $A_{2n+1}$ are the unknowns.
In the case of the periodic Schr\"{o}dinger operators the Novikov--Veselov equations preserve the spectral curve of $L$.

\bigskip
\noindent
{\bf Example 1.}
Let $n=1$. Then
$$
L_3=\partial^3_{z}+u \partial_z +\partial_{\bar{z}}^3+w\partial_{\bar z}, \qquad
B_0=u_z+w_{\bar z},
$$
where
$$ u_{\bar z}=3V_z,\quad\quad\quad w_z=3V_{\bar{z}}.  $$
The first Novikov-Veselov equation has the following form
$$
 V_t=V_{zzz}+V_{\bar{z}\bar{z}\bar{z}}+uV_z+wV_{\bar{z}}+Vw_{\bar{z}}+Vu_z.
$$

\bigskip
\noindent
{\bf Example 2.}
In the case $n=2$, we have
$$
L_5=\partial_z^5+u_3 \partial_z^3+{u_3}_z\partial_z^2+u_1 \partial_z+\partial_{\bar z}^5v_3\partial_{\bar z}^3+{v_3}_{\bar z}\partial_{\bar z}^2+v_1 \partial_{\bar z},
$$
$$
 B_2={u_3}_{z} \partial_z^2+{u_3}_{zz}\partial_z+{w_3}_{\bar z}\partial_{\bar z}^2+{w_3}_{{\bar z}{\bar z}}\partial_{\bar z}+{u_1}_{z}+{w_1}_{\bar z},
$$
where $u_j,w_j$ satisfies the equations.
\begin{equation}\label{eq11}
{u_3}_{\bar z}=5V_z, \qquad
{w_3}_{\bar z}=5V_{\bar z},
\end{equation}
\begin{equation}\label{eq12}
{u_1}_{\bar z}=10 V_{zzz}+3u_3 V_z+{u_3}_{z} V-{u_3}_{zz\bar z},
\end{equation}
\begin{equation}\label{eq13}
{w_1}_{z}=10 V_{\bar{z}\bar{z}\bar{z}}+3w_3 V_{\bar z}+{w_3}_{\bar z}V-{w_3}_{z{\bar z}{\bar z}}.
\end{equation}
The second Novikov--Veselov equation has the following form
$$
 V_t=\partial_z^5 V+u_3 V_{zzz}+2 {u_3}_z V_{zz}+(u_1+{u_3}_{zz}) V_z +
$$
$$
 \partial_{\bar z}^5 V+w_3 V_{{\bar z}{\bar z}{\bar z}}
+2 {w_3}_{\bar z} V_{\bar{z}\bar{z}}+(w_1+{w_3}_{{\bar z}{\bar z}}) V_z +({w_1}_{\bar z}+{u_1}_z)V.
$$
\bigskip

It turns out that if we assume that $V=e^{v(z, {\bar z})}$, where $v$ is a real function satisfying the  Tziz\'eica equation, then the
first stationary Novikov--Veselov equation follows from the  Tziz\'eica equation.

\bigskip
\noindent
{\bf Theorem 4.} \cite{M2}
{\it The stationary Novikov--Veselov equation
$$[A_3, L]+B_0 L=0,$$
with $A_3=\partial_z^3+\partial^3_{\bar z}-(v_z^2+v_{zz})\partial_z-(v_{\bar{z}}^2 +v_{\bar{z}\bar{z}})\partial_{\bar z}$,
$B_0= -\partial_z(v_z^2+v_{zz})-\partial_{\bar z} (v_{\bar{z}}^2 +v_{\bar{z}\bar{z}})$, is equivalent to the following  equations
$$\partial_z(e^{-2v}-e^v-v_{z\bar z})+2v_z(e^{-2v}-e^v-v_{z\bar z})=0,
$$
$$
\partial_{\bar z}(e^{-2v}-e^v-v_{z\bar z})+2v_{\bar z}(e^{-2v}-e^v-v_{z\bar z})=0,
$$
which follows from the Tziz\'eica equation.}

Moreover, coefficients of $A_{2n+1}$ can be expressed in terms of $v$ and derivatives of $v$.
In particular it means that the equations (\ref{eq11})--(\ref{eq13}) can be solved explicitly.

\bigskip
\noindent
{\bf Theorem 5.} \cite{M2}
{\it
Let a real function $v$ satisfies the Tziz\'{e}ica equation and the functions
\begin{eqnarray*}
& & V=e^v,\\
& & u_3=-\frac{5}{3}(v_z^2+v_{zz}), \quad w_3=\bar{u}_3,\\
& & u_1=\frac{5}{9}v_z^4+\frac{10}{9}v_z^2 v_{zz}-\frac{5}{3}v_z^2-\frac{20}{9}v_zv_{zzz}-\frac{10}{9} v_{zzzz},
\quad\ w_1=\bar{u}_1,
\end{eqnarray*}
satisfy equations (\ref{eq11})--(\ref{eq13}).
The second Novicov--Veselov equation attains the following form
$$
 V_t=h+\bar{h},
$$
$$
 h=\frac{1}{9}(5v_1v_2^2+5v_1^2v_3-5v_2v_3-v_1^5),
$$
where $v_j=\partial_z^j v.$
}

\bigskip

If we rewrite the Novikov--Vesleov deformation
$$\partial_{t_2} r=A_5 r$$
in terms of the frame $R$, we get the corresponding matrix $T\in {\rm su}(3)$ such that $R_t=TR$.
The matrix $T$ is given by (\ref{eq:C}), where
$$
a_1=\frac{1}{72\sqrt{2}} e^{-\frac{v}{2}} (144-e^v v_y^4-3e^v v_y^2-e^v v_{yyyy}-e^v v_x^4+12 e^v v_{xy}^2-3e^v v_x v_{xyy})
$$
$$
 -2e^v v_{yy}(2v_x^2-3v_{xxx}+2e^v v_y^2(2v_{yy}+3v_x^2-2v_{xx})+4e^v v_x^2 v_{xx}-3e^v v_{xx}^2+e^v v_{y}(v_{yyy}
$$
$$
 -16v_x v_{xy}-3v_{xxy})+6e^v v_{xxyy}+e^v v_x v_{xxx}-e^v v_{xxxx}),
$$
$$
 a_2=\frac{1}{18\sqrt{2}}e^{-\frac{3v}{2}}(-e^{2v}v_y^3 v_x+2e^{2v}v_y^2v_{xy}+2v_{xy}(-2+8e^{3v}-e^{2v}v_x^2+3e^{2v}v_{xx}))
$$
$$
 +v_y(e^{2v}v_x^3-4v_x(4+2e^{3v}+e^{2v}v_{xx})-e^{2v}v_{xxx})+e^{2v}(-v_x v_{xxy}+2v_{xxxy}),
$$
$$
 s=\frac{1}{3}(v_x^2-v_y^2+v_{xx}-v_{yy}).
$$
Previous formula for $s$ corresponds to an eigenfunction $s_1$ of the Laplace--Beltrami operator with eigenvalue $6$ in Theorem 3.
By direct calculation one can check that $s_2$ is also an eigenfunction with the same eigenvalue and
functions $a_1, a_2$ and $s$ give a solution of the equations \eqref{eq5}--\eqref{eq10}.



\end{document}